\documentclass[seceqn,dvips]{arxbj}

\aid{0}
\volume{16}
\issue{4}
\pubyear{2010}
\firstpage{1369}
\lastpage{1384}
\doi{10.3150/10-BEJ252}

\makeatletter

\newtheorem{theorem}{Theorem}[section]

\newtheorem{corollary}{Corollary}[section]

\newremark{example}{Example}
\newremark{remark}{Remark}[section]

\def\eqref#1{(\ref{#1})}

\makeatother

\begin{document}
\begin{frontmatter}

\title{Consistent group selection in high-dimensional linear regression}
\runtitle{Consistent group selection}

\begin{aug}
\author[1]{\fnms{Fengrong} \snm{Wei}\thanksref{1}\ead[label=e1]{fwei@westga.edu}\corref{}} \and
\author[2]{\fnms{Jian} \snm{Huang}\thanksref{2}\ead[label=e2]{jian-huang@uiowa.edu}}
\runauthor{F. Wei and J. Huang}
\address[1]{Department of Mathematics, University of West Georgia,
1601 Maple Street, Carrollton, GA 30118, USA. \printead{e1}}
\address[2]{Department of Statistics and Actuarial Science, University
of Iowa, Iowa City, IA 52242, USA.\\ \printead{e2}}
\end{aug}

% HISTORY:
\received{\smonth{11} \syear{2008}}
\revised{\smonth{9} \syear{2009}}

% ABSTRACT
%
\begin{abstract}
In regression problems where covariates can be naturally grouped,
the group Lasso is an attractive method for variable selection since it
respects the grouping structure in the data.
We study the selection and estimation properties of the group Lasso in
high-dimensional
settings when the number of groups exceeds the sample size. We provide
sufficient conditions
under which the group Lasso selects a model whose dimension is
comparable with the underlying
model with high probability and is estimation consistent. However, the
group Lasso is, in general,
not selection consistent and also tends to select groups that are not
important in the model.
To improve the selection results, we propose an adaptive group Lasso
method which is a
generalization of the adaptive Lasso and requires an initial estimator.
We show that the adaptive group Lasso is consistent in group selection
under certain conditions
if the group Lasso is used as the initial estimator.
\end{abstract}

% KEYWORDS
%

\end{frontmatter}

%s1 ###
\section{Introduction}

Consider the linear regression model with $p$ groups of covariates
\begin{eqnarray*}
Y_{i}=\sum_{k=1}^{p}X_{ik}'\beta_{k}+\varepsilon_{i} , \qquad i=1,\ldots
,n ,
\end{eqnarray*}
where $Y_{i}$ is the response variable, $\varepsilon_{i}$ is the
error term, $X_{ik}$ is a $d_k \times1$ covariate vector
representing the $k$th group and $\beta_{k}$ is the corresponding
$d_k \times1$ vector of regression coefficients. For such a model, the
group Lasso (Antoniadis and Fan (\citeyear{AF2001}), Yuan and Lin (\citeyear{YL2006})) is an
attractive method for variable selection since it respects the
grouping structure in the covariates. This method is a natural
extension of the Lasso (Tibshirani (\citeyear{T1996})), in which an $\ell_2$-norm
of the coefficients associated with a group of variables is used as
a component in the penalty function. However, the group Lasso is, in
general, not selection consistent and tends to select more groups
than there are in the model. To improve the selection results, we
consider an adaptive group Lasso method which is a generalization
of the adaptive Lasso (Zou (\citeyear{Z2006})). We provide sufficient conditions
under which the adaptive group Lasso is selection consistent if the
group Lasso is used as the initial estimator.

The need to select groups of variables arises in many statistical
modeling problems and applications. For example, in multifactor analysis
of variance, a factor with multiple levels can be represented by a
group of dummy variables. In nonparametric additive regression, each
component can be expressed as a linear combination of a set of basis
functions. In both cases, the selection of important factors or
nonparametric components amounts to the selection of groups of
variables. Several recent papers have considered group selection
using penalized methods. In addition to the group Lasso, Yuan and
Lin (\citeyear{YL2006}) have proposed the group Lars and group non-negative garrote
methods. Kim, Kim and Kim (\citeyear{KKK2006}) considered the group Lasso in the
context of generalized linear models. Zhao, Rocha and Yu (\citeyear{ZRY2008})
proposed a composite absolute penalty for group selection, which can
be considered a generalization of the group Lasso. Meier, van de
Geer and B\"{u}hlmann (\citeyear{MvB2008}) studied the group Lasso for logistic
regression. Huang, Ma, Xie and Zhang (\citeyear{HHM2008}) proposed a group bridge
method that can be used for simultaneous group and individual
variable selection.

There has been much work on the penalized methods for variable
selection and estimation with high-dimensional data. Several
approaches have been proposed, including the least absolute
shrinkage and selection operator (Lasso, Tibshirani (\citeyear{T1996})), the
smoothly clipped absolute deviation (SCAD) penalty (Fan and Li (\citeyear{FL2001}),
Fan and Peng (\citeyear{FP2004})), the elastic net (Enet) penalty (Zou and Hastie
(\citeyear{ZH2006})) and the minimum concave penalty (Zhang (\citeyear{Z2007})). Much progress
has been made in understanding the statistical properties of these
methods in both fixed $p$ and $p \gg n$ settings. In particular,
several recent studies considered the Lasso with regard to its
variable selection, estimation and prediction properties; see, for
example, Knight and Fu (\citeyear{KF2001}), Greenshtein and Ritov (\citeyear{GR2004}),
Meinshausen and Buhlmann (\citeyear{MB2006}), Zhao and Yu (\citeyear{ZY2006}), Huang, Ma and
Zhang (\citeyear{HMZ2006}), van de Geer (\citeyear{v2008}) and Zhang and Huang (\citeyear{ZH2008}), among
others. All of these studies are concerned with the Lasso for
individual variable selection.

In this article, we study the asymptotic properties of the group Lasso
and the adaptive group Lasso in high-dimensional settings when $p \gg n$.
We generalize the results concerning the Lasso obtained in Zhang and
Huang (\citeyear{ZH2008}) to the group Lasso.
We show that, under a generalized sparsity condition and the sparse
Riesz condition, as well as certain regularity conditions, the group
Lasso selects a model whose dimension has the same order as the
underlying model,
%controls the bias of the selected model at a level determined by the
%%contributions of small regression coefficients and threshold bias, and
selects all groups whose $\ell_2$-norms are of greater order than
the bias of the selected model and is estimation consistent. In
addition, under a narrow-sense sparsity condition (see page \pageref{p4}) and
using the
group Lasso as the initial estimator, the adaptive group Lasso can
correctly select important groups with high probability.

Our theoretical and simulation results suggest the following
one-step approach to group selection in high-dimensional settings.
First, we use the group Lasso to obtain an initial estimator and
reduce the dimension of the problem. We then use the adaptive group
Lasso to select the final set of groups of variables. Since the
computation of the adaptive group Lasso estimator can be carried out
using the same algorithm and program for the group Lasso, the
computational cost of this one-step approach is approximately twice
that of a single group Lasso computation. This approach, iteratively
using the group Lasso twice, follows the idea of the adaptive Lasso
(Zou (\citeyear{Z2006})) and a proposal by B\"{u}hlmann and Meier (\citeyear{BM2008}) in the
context of individual variable selection.

The rest of the paper is organized as follows. In Section~\ref{sec2}, we state
the results on the selection, bias of the selected model and convergent
rate of the group Lasso estimator. In Section~\ref{sec3}, we describe the
selection and estimation consistency results concerning the adaptive
group Lasso. In Section~\ref{sec4}, we use simulation to compare the group
Lasso and adaptive group Lasso. Proofs are given in Section~\ref{sec5}.
Concluding remarks are given in Section~\ref{sec6}.

%s2 ###
\section{The asymptotic properties of the group Lasso}\label{sec2}

Let $Y=(Y_1, \ldots, Y_n)'$ and $X=(X_{1},\ldots,X_{p})$, where
$X_{k}$ is the $n \times d_{k}$ covariate
submatrix corresponding to the $k$th group. For a given penalty level
$\lambda\geq0$, the group Lasso
estimator of $\beta= (\beta_1', \ldots, \beta_p')'$ is
%
%e2.1 ###
\begin{equation}\label{eq:GL1}
\hat{\beta}=\arg\min_{\beta}\frac{1}{2}(Y-X\beta)^{\mathrm{T}}(Y-X\beta
)+\lambda\sum_{k=1}^{p}\sqrt{d_{k}}\|\beta_{k}\|_{2} ,
\end{equation}
where ${\hat{\beta}}=(\hat{\beta}_{1}^{\prime},\ldots,\hat{\beta
}_{p}^{\prime
})^{\prime}$.

We consider the model selection and estimation properties of
$\hat{\beta}$ under a generalized sparsity condition (GSC) of the model
and a sparse Riesz condition (SRC) on the covariate matrix. These
two conditions were first formulated in the study of the Lasso
estimator (Zhang and Huang (\citeyear{ZH2008})). The GSC assumes that for some
$\eta_1 \ge0$, there exists an $A_{0} \subset
\{1,\ldots, p\}$ such that $\sum_{k \in
A_{0}}\|\beta_{k}\|_{2} \le\eta_{1}$, where $\|\cdot\|_2$ denotes the
$\ell_2$-norm. Without loss of
generality, let $A_{0}=\{q+1,\ldots,p\}$. The GSC is
then
%
%e2.2 ###
\begin{equation}
\label{GSCa} \sum_{k=q+1}^{p}\|\beta_{k}\|_{2}\le\eta_{1} .
\end{equation}
The number of truly important groups is thus $q$. A more rigid
way to describe sparsity is to assume $\eta_1=0$, that is,\label{p4}
%
%e2.3 ###
\begin{equation}
\label{NSCa} \|\beta_{k}\|_{2}=0 ,\qquad   k=q+1, \ldots, p .
\end{equation}
This is a special case of the GSC and we call it the \emph
{narrow-sense sparsity condition} (NSC). In practice, the GSC is a more
realistic formulation of a sparse model. However, the NSC can often be
considered
a reasonable approximation to the GSC, especially when $\eta_1$ is smaller
than the noise level associated with model fitting.

The SRC controls the range of eigenvalues of the submatrix. For $A
\subset\{1, \ldots, p\}$, we define $X_{A}=(X_{k}, k \in A)$ and
$\Sigma_{AA}=X_{A}^{\prime}X_{A}/n$. Note that $X_{A}$ is an $n
\times\sum_{k\in A} d_k$ matrix. The design matrix $X_{A}$
satisfies the sparse Riesz condition (SRC) with rank $q^{*}$ and
spectrum bounds $0 < c_{*}<c^{*} < \infty$ if
%
%e2.4 ###
\begin{equation}
\label{SRCa}
c_{*} \leq\frac{\|X_{A}\nu\|_{2}^{2}}{n\|\nu\|_{2}^{2}}
\leq c^{*}\qquad \forall A \mbox{ with }q^{*}=|A|=\#\{k\dvt k \in A\}
\mbox{ and } \nu\in R^{\sum_{k \in A}d_{k}} .
\end{equation}
%
%where $|A|=\sum_{k \in A}d_{k}$.

Let ${\hat{A}}=\{k\dvt \|\hat{\beta}_{k}\|_{2}>0, 1 \le k \le p\}$,
which is the
set of indices of the groups selected by the group Lasso. An
important quantity is the cardinality of $\hat{A}$, defined as
%
%e2.5 ###
\begin{equation}
\label{q}
\hat{q}=|\hat{A}|=\#\{k\dvt \|\hat{\beta}_{k}\|_{2}>0, 1 \le k \le p\},
\end{equation}
which determines the dimension of the selected model. If
$\hat{q}=\mathrm{O}(q)$, then the selected model has dimension comparable to the
underlying model. Following Zhang and Huang (\citeyear{ZH2008}), we also consider
two measures of the selected model. The first measures the error of
the selected model:
%
%e2.6 ###
\begin{equation}\label{w}
\tilde{\omega}=\|(I-{\hat{P}})X\beta\|_{2} ,
\end{equation}
where ${\hat{P}}$ is the projection matrix from $R^{n}$ to the linear span
of the set of selected groups and $I \equiv I_{n \times n}$ is the
identity matrix. Thus, $\tilde{\omega}^{2}$ is the sum of squares of
the mean vector not accounted for by the selected model. To measure the
important groups missing in the selected model, we define
%
%e2.7 ###
\begin{equation}\label{zeta}
\zeta_{2}=\biggl(\sum_{k \notin
A_{0}}\|\beta_{k}\|_{2}^{2}I\{\|\hat{\beta}_{k}\|_{2}=0\}\biggr)^{1/2} .
\end{equation}

We now describe several quantities that will be useful in describing
the main results. Let $d_{a}=\max_{1\le k\le p}d_{k}$,
$d_{b}=\min_{1\le k\le p}d_{k}$, $d=d_{a}/d_{b}$ and
$N_{d}=\sum_{k=1}^{p}d_{k}$. Define
%
%e2.8 ###
\begin{equation}\label{eq:r1r2GL}
r_{1}\equiv r_{1}(\lambda)=\biggl(\frac{nc^{*}\sqrt{d_{a}}\eta
_{1}}{\lambda d_{b}q}\biggr)^{1/2} , \qquad r_{2}\equiv
r_{2}(\lambda) =\biggl(\frac{nc^{*}\eta_{2}^{2}}{\lambda
^{2}d_{b}q}\biggr)^{1/2} , \qquad\bar{c}=\frac{c^{*}}{c_{*}} ,
\end{equation}
where $\eta_{2}\equiv\max_{A \subset A_{0}}\|\sum_{k \in
A}X_{k}\beta_{k}\|_{2}$,
%
%e2.11 ###
%e2.10 ###
%e2.9 ###
\begin{eqnarray}
\label{eq:M1GL}
M_{1}&\equiv& M_{1}(\lambda)=2+4r_{1}^{2}+4\sqrt{d\bar{c}}
r_2+4d\bar{c} ,\\
\label{eq:M2GL}
M_{2}&\equiv& M_{2}(\lambda)=\tfrac{2}{3}\bigl(1+4r_{1}^{2}+2d\bar{c}+
4\sqrt{2d}\bigl(1+\sqrt{\bar{c}}\bigr) \sqrt{\bar{c}}  r_2
+\tfrac{16}{3}d\bar{c}^{2}\bigr) ,
\\
\label{eq:M3GL}
M_{3}&\equiv&
M_{3}(\lambda)=\tfrac{2}{3}\bigl(1+4r_{1}^{2}+
4\sqrt{d\bar{c}}\bigl(1+2\sqrt{1+\bar{c}} \bigr)r_{2}+3r_{2}^{2}
+\tfrac{2}{3}d\bar{c}(7+4\bar{c})\bigr) .
\end{eqnarray}

Let
$\lambda_{n,p}=2\sigma\sqrt{8(1+c_{0})d_{a}d^{2}q^{*}\bar
{c}nc^{*}\log(N_{d}\vee
a_{n})}$, where $c_{0}\geq0$ and $a_{n}\geq0$, satisfying
$pd_{a}/(N_{d}\vee a_{n})^{1+c_{0}}\approx0$, and $\lambda_{0}=\inf
\{\lambda\dvt M_{1}q+1 \leq q^{*}\}$, where $\inf\varnothing= \infty$. We
also consider the constraint
%
%e2.12 ###
\begin{equation}\label{eq:Crk2GL}
\lambda\geq\max\{\lambda_{0},\lambda_{n,p}\} .
\end{equation}
For large $p$, the lower bound here is allowed to be
$\lambda_{n,p}=2\sigma[8(1+c_{0})d_{a}d^{2}q^{*}\bar{c}nc^{*}\log
(N_{d})]^{1/2}$ with $a_{n}=0$; for fixed $p$,
$a_{n}\rightarrow\infty$ is required.

We assume the following basic condition.
\begin{enumerate}
\item[(C1)]
The errors $\varepsilon_{1},\ldots, \varepsilon_{n}$ are independent
and identically distributed as $N(0, \sigma^2)$.
\end{enumerate}

\begin{theorem}
\label{ThmA}
Suppose that $q \geq1$ and that \textup{(C1)}, the GSC \textup{(\ref{GSCa})} and SRC
\textup{(\ref{SRCa})} are satisfied. Let $\hat{q},\tilde{\omega}$ and
$\zeta_{2}$ be defined as in $\eqref{q}$, $\eqref{w}$ and
$\eqref{zeta}$, respectively, for the model $\hat{A}$ selected by the
group Lasso from $\eqref{eq:GL1}$. Let $M_{1},M_{2}$ and $M_{3}$ be
defined as in $\eqref{eq:M1GL}$, $\eqref{eq:M2GL}$ and
$\eqref{eq:M3GL}$, respectively. If the constraint
$\eqref{eq:Crk2GL}$ is satisfied, then the following assertions hold with
probability converging to \textup{1}:
\begin{eqnarray*}
{\hat{q}} & \leq&\#\{k\dvt \|\hat{\beta}_{k}\|_{2} >0   \mbox{ or }
  k \notin
A_{0}\} \leq M_{1}(\lambda)q ,
\\
\tilde{\omega}^{2} &=&\|(I-\hat{P})X\beta\|_{2}^{2} \leq
M_{2}(\lambda)B_{1}^{2}(\lambda) ,
\\
\zeta_{2}^{2} & =&\sum_{k \notin
A_{0}}\|\beta_{k}\|_{2}^{2}I\{\|\hat{\beta}_{k}\|_{2}=0\}\leq\frac{
M_{3}(\lambda)B_{1}^{2}(\lambda)}{c_{*}n} ,
\end{eqnarray*}
where $B_{1}(\lambda)=((\lambda^{2}d_{b}^{2}q)/(nc^{*}))^{1/2}$ .
\end{theorem}

\begin{remark}
\label{Remark0} The condition $q \geq1$ is not necessary since it
is only used to express quantities in terms of ratios in
\textup{(\ref{eq:r1r2GL})} and Theorem \ref{ThmA}. If $q=0$, we use $
r_{1}^{2}q = nc^{*}\sqrt{d_{a}}\eta_{1}/(\lambda d_{b})$ and
$r_{2}^{2}q=nc^{*}\eta_{2}^{2}/(\lambda^{2}d_{b})$ to recover
$M_{1}$, $M_{2}$ and $M_{3}$ in (\ref{eq:M1GL}), (\ref{eq:M2GL}),
(\ref{eq:M3GL}), respectively, giving the results $\hat{q}\leq
4nc^{*}\sqrt{d_{a}}\eta_{1}/\lambda d_{b}$, $\tilde{\omega}^{2}
\leq
8\lambda\sqrt{d_{a}}d_{b}\eta_{1}/3$ and $\zeta_{2}^{2}=0$.
\end{remark}

\begin{remark}
\label{Remark1}
If $\eta_1=0$ in (\ref{GSCa}), then $r_{1}=r_{2}=0$ and
\begin{eqnarray*}
M_{1}=2+4d\bar{c}, \qquad  M_{2}=\tfrac{2}{3}\bigl(1+2d\bar{c}
+\tfrac{16}{3}d\bar{c}^{2}\bigr), \qquad  M_{3}=\tfrac{2}{3}\bigl(1
+\tfrac{2}{3}d\bar{c}(7+4\bar{c})\bigr),
\end{eqnarray*}
all of which depend only on $d$ and $\bar{c}$. This suggests that the
relative sizes of the groups affect the selection
results. Since $d \ge1$, the most favorable case is $d=1$, that is,
when the groups have equal sizes.
\end{remark}

\begin{remark}
\label{Remark2} If $d_1 = \cdots= d_p=1$, the group Lasso
simplifies to the Lasso and Theorem \ref{ThmA} is a direct
generalization of Theorem 1 on the selection properties of the Lasso
obtained by Zhang and Huang (\citeyear{ZH2008}). In particular, when $d_1 =
\cdots= d_p=1$, $r_{1}, r_{2}, M_{1}, M_{2}, M_{3}$
are the same as the constants in Theorem 1 of Zhang and Huang (\citeyear{ZH2008}).
\end{remark}

\begin{remark}
\label{Remark3}
A more general definition of the group Lasso is
%
%e2.13 ###
\begin{equation}\label{eq:GL0}
\hat{\beta}{}^{*}=\arg\min_{\beta}\frac{1}{2}(Y-X\beta
)^{\prime}(Y-X\beta)+\lambda
\sum_{k=1}^{p}(\beta_{k}^{\prime}R_{k}\beta_{k})^{1/2} , %\|
\end{equation}
where $R_{k}$ is a $d_{k} \times d_{k}$ positive definite matrix.
This is useful when certain relationships among the coefficients
need be specified. By the Cholesky decomposition, there exists a
matrix $Q_k$ such that $R_{k}=d_{k}Q_{k}^{\prime}Q_{k}$. Let
$\beta^{*}=Q_{k}\beta$, and $X_{k}^{*}=X_{k}Q_{k}^{-1}$. Then,
(\ref{eq:GL0}) becomes
\begin{eqnarray*}
\hat{\beta}{}^{*}=\arg\min_{\beta^{*}}(Y-X^{*}\beta
^{*})^{\prime}(Y-X^{*}\beta^{*})+\lambda
\sum_{k=1}^{p}\sqrt{d_{k}}\|\beta_{k}^{*}\|_{2} .
\end{eqnarray*}
The GSC for (\ref{eq:GL0}) is $\sum_{k=q+1}^{p}
(\beta_{k}^{\prime}Q_{k}^{\prime}Q_{k}\beta_{k})^{1/2} \le\eta_1$. The
SRC can be assumed for $X \cdot Q^{-1}$, where $X \cdot
Q^{-1}=(X_{1}Q_{1}^{-1},\ldots, X_{p}Q_{p}^{-1})$.
\end{remark}

Immediately, from Theorem \ref{ThmA}, we have the following corollary.
\begin{corollary}
\label{ColA}
Suppose that the conditions of Theorem \ref{ThmA} hold and $\lambda$
satisfies the constraint (\ref{eq:Crk2GL}).
Then, with probability converging to one,
all groups with $\|\beta_k\|_2^2 > M_3(\lambda) q \lambda
^2/(c_{*}c^{*}n^2)$ are selected.
\end{corollary}

From Theorem \ref{ThmA} and Corollary \ref{ColA}, the group Lasso
possesses similar properties to the Lasso in terms of sparsity and bias
(Zhang and Huang (\citeyear{ZH2008})). In particular, the group Lasso selects a model
whose dimension has the same order as the underlying model.
Furthermore, all of the groups with coefficients whose $\ell_2$-norms
are greater than the threshold given in Corollary \ref{ColA} are
selected with high probability.

\begin{theorem}
\label{ThmB} Let $\{\bar{c},\sigma, r_{1}, r_{2}, c_{0}, d\}$ be
fixed and $1 \leq q \leq n \leq p\rightarrow\infty$. Suppose that
the conditions in Theorem \ref{ThmA} hold. Then, with probability
converging to 1,
we have
\begin{eqnarray*}
\|\hat{\beta}-\beta\|_{2} \leq\frac{1}{\sqrt{nc_{*}}}
\bigl(2\sigma
\sqrt{M_{1}\log(N_{d})q}+
\bigl(r_{2}+\sqrt{dM_{1}\bar{c}}\bigr)B_{1}\bigr)+\sqrt{\frac
{c_{*}r_{1}^{2}+r_{2}^{2}}{c_{*}c^{*}}}\frac{\sqrt{q}\lambda}{n}
\end{eqnarray*}
and
\begin{eqnarray*}
\|X\hat{\beta}-X\beta\|_{2} \leq2\sigma\sqrt{M_{1}\log
(N_{d})q}+\bigl(2r_{2}+ \sqrt{dM_{1}\bar{c}}\bigr)B_{1} .
\end{eqnarray*}
\end{theorem}

Theorem \ref{ThmB} is stated for a general $\lambda$ that satisfies
(\ref{eq:Crk2GL}). The following result is an immediate corollary of
Theorem \ref{ThmB}.

\begin{corollary}
\label{ColB} Let $\lambda= 2\sigma\sqrt
{8(1+c_{0}^{\prime})d_{a}d^2q^{*}\bar{c}c^{*}n\log(N_{d})}$ with a fixed
$c_{0}^{\prime} \geq c_{0}$. Suppose that all of the conditions in Theorem
\ref{ThmB} hold. We then have
\[
\|\hat{\beta}-\beta\|_{2} = \mathrm{O}_p\bigl(\sqrt{q\log(N_{d})/n} \bigr)
\quad
\mbox{and}\quad
\|X\hat{\beta}-X\beta\|_{2}=\mathrm{O}_p\bigl(\sqrt{q\log(N_{d})} \bigr).
\]
\end{corollary}

This corollary follows by substituting the given $\lambda$ value
into the expressions in the results of Theorem \ref{ThmB}.

%s3 ###
\section{Selection consistency of the adaptive group Lasso}\label{sec3}

As shown in the previous section, the group Lasso has excellent
selection and estimation properties. However, there is room for
improvement, particularly with regard to selection. Although the group
Lasso selects a model whose dimension is comparable to that of
the underlying model, the simulation results reported in Yuan and Lin
(\citeyear{YL2006}) and those reported below suggest that it tends to select more
groups than there are in the underlying model. To correct the
tendency of overselection by the group Lasso, we generalize the
idea of the adaptive Lasso (Zou (\citeyear{Z2006})) for individual variable
selection to the present problem of group selection.

Consider a general group Lasso criterion with a weighted penalty term,
%
%e3.1 ###
\begin{equation}\label{eq:ModelAGL}
\frac{1}{2}(Y-X\beta)^{\prime}(Y-X\beta) +\tilde{\lambda} \sum
_{k=1}^{p} w_{k}\sqrt{d_{k}}\|\beta_{k}\|_{2} ,
\end{equation}
where $w_{k}$ is the weight associated with the $k$th group. The
$\lambda_k \equiv\tilde{\lambda} w_{k}$ can
be regarded as the penalty level corresponding to the $k$th group. For
different groups,
the penalty level $\lambda_{k}$ can be different. If
we can have lower penalty for groups with large coefficients and higher
penalty for groups with small
coefficients (in the $\ell_2$ sense), then we expect to be able to
improve variable selection accuracy and reduce
estimation bias. One way to obtain the information about whether a
group has large or small coefficients is by using a consistent initial
estimator.

Suppose that an initial estimate ${\tilde{\beta}}$ is available.
A simple approach to determining the weight is to use the initial
estimator. Consider
%
%e3.2 ###
\begin{equation}
\label{AdaW}
w_{k}=\frac{1}{\|{\tilde{\beta}_{k}}\|_2} ,\qquad  k=1, \ldots, p .
\end{equation}
Thus, for each group, its penalty is proportional to the inverse of
the norm of $\tilde{\beta}_k$. This choice of the penalty level for each
group is a natural generalization of the adaptive Lasso (Zou (\citeyear{Z2006})).
In particular, when each group only contains a single variable,
(\ref{AdaW}) simplifies to the adaptive Lasso penalty.

Let $\theta_{a}=\max_{k \in A_{0}^{c}}\|\beta_{k}\|_{2}$ and $\theta
_{b}=\min_{k \in A_{0}^{c}}\|\beta_{k}\|_{2}$.
We say that an initial estimator $\tilde{\beta}$ is consistent at
zero with
rate $r_{n}$ if $r_{n}\max_{k \in A_{0}}\|\tilde{\beta}_{k}\|_{2}=\mathrm{O}_{p}(1)$,
where $r_{n} \rightarrow
\infty$ as $n \rightarrow\infty$, and there exists a constant $\xi
_{b}>0$ such that for any $\varepsilon>0$, $P(\min_{k \in
A_{0}^{c}}\|\tilde{\beta}_{k}\|_{2} > \xi_{b}\theta_{b})>
1-\varepsilon$
for $n$ sufficiently large.

In addition to (C1), we assume the following conditions:
\begin{enumerate}
\item[(C2)]
the initial estimator $\tilde{\beta}$ is consistent at zero with rate $r_{n}
\rightarrow\infty$;
\item[(C3)]
\begin{eqnarray*}\\[-25pt]
\frac{\sqrt{d_{a}(\log q)}}{\sqrt{n}\theta_{b}}\rightarrow0 ,
\qquad
\frac{\tilde{\lambda} d_{a}^{3/2}q}{n\theta_{b}^{2}}\rightarrow0
,\qquad
\frac{\sqrt{nd\log(p-q)}}{\tilde{\lambda} r_{n}} \rightarrow0 ,
\qquad
\frac{d_{a}^{5/2}q^{2}}{r_{n}\theta_{b}\sqrt{d_{b}}}\rightarrow0 ;
\end{eqnarray*}
\item[(C4)]
all of the eigenvalues of $\Sigma_{A_{0}^{c}A_{0}^{c}}$ are bounded
away from zero and infinity.
\end{enumerate}

Condition (C2) assumes that an initial zero-consistent estimator
exists. It is the most critical one and is generally difficult to
establish. It assumes that we can consistently differentiate between
important and non-important groups. For fixed $p$ and $d_k$, the
ordinary least-squares estimator can be used as the initial
estimator. However, when $p > n$, the least-squares estimator is no
longer feasible. By Theorems \ref{ThmA} and \ref{ThmB}, the group
Lasso estimator $\hat{\beta}$ is consistent at zero with rate
$\sqrt{n/(q\log(N_{d}))}$. Condition (C3) restricts the numbers of
important and non-important groups, as well as variables within the
groups. It also places constraints on the penalty parameter and the
$\ell_{2}$-norm of the smallest important group. Condition (C4)
assumes that the eigenvalues of $\Sigma_{A_{0}^{c}A_{0}^{c}}$ are
finite and bounded away from zero. This is reasonable since the
number of important groups is small in a sparse model. This
condition ensures that the true model is identifiable.

Define
%
%e3.3 ###
\begin{equation}\label{eq:ModelAGLb}
\hat{\beta}^*=\arg\min\frac{1}{2}(Y-X\beta)^{\prime}(Y-X\beta)
+\tilde
{\lambda} \sum_{k=1}^{p} \|\tilde{\beta}_k\|_2^{-1}
\sqrt{d_{k}}\|\beta_{k}\|_{2} .
\end{equation}

\begin{theorem}
\label{ThmC} If \textup{(C1)--(C4)} and NSC \textup{(\ref{NSCa})} are
satisfied, then
\begin{eqnarray*}
P(\|\hat{\beta}_{k}^{*}\|_{2} \neq0 ,k \notin A_0, \|\hat
{\beta}
_{k}^{*}\|_{2}=0, k \in A_0)\rightarrow1 .
\end{eqnarray*}
\end{theorem}

Therefore, the adaptive group Lasso is selection consistent if the
conditions stated in Theorem~\ref{ThmA} hold.

If we use $\hat{\beta}$ as the initial estimator, then (C3) can be
changed to
\begin{enumerate}
\item[(C3)$^*$]
\begin{eqnarray*}\\[-25pt]
\frac{\sqrt{d_{a}(\log q)}}{\sqrt{n}\theta_{b}} &\rightarrow& 0 ,\qquad
\frac{\tilde{\lambda} d_{a}^{3/2}q}{n\theta_{b}^{2}}\rightarrow0 ,\qquad
  \frac{\sqrt{dq\log(p-q)\log(N_{d})}}{\tilde{\lambda}}
\rightarrow0 ,
\\
\frac{(d_{a}q)^{5/2}\sqrt{\log(N_{d})}}{\theta_{b}\sqrt
{nd_{b}}}&\rightarrow&0 .
\end{eqnarray*}
\end{enumerate}
We often have $\tilde{\lambda}=n^\alpha$ for some $0<\alpha<1/2$. In
this case, the number of non-important groups can be as large as
$\exp(n^{2\alpha}/(q\log q))$ with the number of important groups
satisfying $q^5\log q/n \rightarrow0$, assuming that $\theta_{b}$
and the number of variables within the groups are finite.

\begin{corollary}
\label{ColC} Let the initial estimator $\tilde{\beta}=\hat{\beta
}$, where
$\hat{\beta}$ is
the group Lasso estimator. Suppose that the NSC \textup{(\ref{NSCa})} holds
and that (\textup{C1), (C2), (C3)$^*$} and \textup{(C4)} are satisfied. We then have
\begin{eqnarray*}
P(\|\hat{\beta}_{k}^{*}\|_{2} \neq0 ,k \notin A_0, \|\hat
{\beta}
_{k}^{*}\|_{2}=0, k \in A_0)\rightarrow1 .
\end{eqnarray*}
\end{corollary}

This corollary follows directly from Theorem \ref{ThmC}. It shows
that the iterated group Lasso procedure that uses a combination of
the group Lasso and the adaptive group Lasso is selection
consistent.

\begin{theorem}
\label{ThmD}
Suppose that the conditions in Theorem \ref{ThmB} hold and that
$\theta_{b}> t_{b}$ for some constant $t_{b}>0$. If $\tilde{\lambda}
\sim \mathrm{O}(n^{\alpha})$ for some $0 < \alpha< 1/2$, then
\begin{eqnarray*}
\|\hat{\beta}^{*}-\beta\|_{2} =
\mathrm{O}_p\Biggl(\sqrt{\frac{q}{n}+\frac{\tilde{\lambda}^{2}}{n^{2}}}\Biggr) =
\mathrm{O}_p\Biggl(\sqrt{\frac{q}{n}}\Biggr) , \qquad\|X\hat{\beta}^{*}-X\beta\|_{2}
\sim
\mathrm{O}\Biggl(\sqrt{q+\frac{\tilde{\lambda}^{2}}{n}}\Biggr) =\mathrm{O}_p\bigl(\sqrt{q}\bigr) .
\end{eqnarray*}
\end{theorem}

Theorem \ref{ThmD} implies that for the adaptive group Lasso, given a
zero-consistent initial estimator, we can
reduce a high-dimensional problem to a lower-dimensional one. The
convergence rate is improved, compared with that of the group Lasso, by
choosing an appropriate penalty parameter~$\tilde{\lambda}$.

%s4 ###
\section{Simulation studies}\label{sec4}

In this section, we use
simulation to evaluate the finite sample performance of the group Lasso
and the adaptive group Lasso. Let $\lambda_{k}=\tilde{\lambda}/\|
\hat{\beta}_{k}\|_{2}$, if
$\|\hat{\beta}_{k}\|_{2}>0$; if
$\|\hat{\beta}_{k}\|_{2}=0$, then $\lambda_{k}=\infty$, $\hat
{\beta}^{*}_{k}=0$.
We can thus drop the corresponding covariates $X_k$ from the
model and only consider the groups with $\|\hat{\beta}_k^{*}\|_2 > 0$.
After a scale transformation, we can directly
apply the group least angle regression algorithm (Yuan and Lin (\citeyear{YL2006})) to
compute the adaptive group Lasso estimator $\hat{\beta}^{*}$.
The penalty parameters for the group Lasso and the adaptive group Lasso
are selected using the BIC criterion
(Schwarz (\citeyear{S1978})).

We consider two scenarios of simulation models. In the first scenario,
the group sizes are equal; in the second, the group sizes vary.
For every scenario, we consider the cases $p < n $ and $p > n$. In
all of the examples, the sample size is $n=200$.

\begin{example}\label{e1} In this example,
there are 10 groups, each consisting of 5 covariates. The
covariate vector is $X=(X_{1},\ldots,X_{10})$, where
$X_{j}=(X_{5(j-1)+1},\ldots,X_{5(j-1)+5})$, $1 \le j \le10$. To
generate~$X$, we first simulate 50 random variables,
$R_{1},\ldots,R_{50}$, independently from $N(0,1)$. Then, $Z_{j}$,
$j=1,\ldots, 10$, are simulated from a
multivariate normal distribution with with mean zero and
$\operatorname{cov}(Z_{j_{1}},Z_{j_{2}})=0.6^{|j_{1}-j_{2}|}$. The covariates
$X_{1},\ldots,X_{50}$ are generated as\looseness=1
\begin{eqnarray*}
X_{5(j-1)+k}=\frac{Z_{j}+R_{5(j-1)+k}}{\sqrt{2}} ,\qquad  1 \leq
j \leq10 ,  1\leq k \leq
5 .
\end{eqnarray*}
The random error $\varepsilon\sim N(0, 3^2)$. The response variable
$Y$ is generated
from $Y=\sum_{k=1}^{10}X_{k}'\beta_{k}+\varepsilon$, where
$\beta_{1}=(0.5,1,1.5,2,2.5),  \beta_{2}=(2,2,2,2,2),   \beta
_{3}=\cdots=\beta_{10}=(0,0,0,0,0)$.
\end{example}

\begin{example}\label{e2} In this example, the number of groups is
$p=10$. Each group consists of 5
covariates. The covariates are generated the same way as in Example
\ref{e1}. However, the regression coefficients
$\beta_{1} =(0.5,1,1.5,1,0.5),   \beta_{2} =(1,1,1,1,1),
\beta_{3} =(-1,0,1,2,1.5),   \beta_{4}= (-1.5,1,0.5,0.5,0.5),
\beta_{5} =\cdots=\beta_{10}=(0,0,0,0,0)$.
\end{example}

\begin{example}\label{e3} In this example, the number of groups
$p=210$ is bigger than the sample size $n$. Each group consists of 5
covariates. The covariates are generated the same way as in Example
\ref{e1}. However, the regression coefficients
$\beta_{1} =(0.5,1,1.5,1,0.5),   \beta_{2} =(1,1,1,1,1),
\beta_{3} =(-1,0,1,2,1.5),   \beta_{4}= (-1.5,1,0.5,0.5,0.5),
\beta_{5} =\cdots=\beta_{210}=(0,0,0,0,0)$.
\end{example}

\begin{example}\label{e4} In this example, the group sizes differ
across groups. There are 5 groups with size 5 and
5 groups with size 3. The covariate vector is
$X=(X_{1},\ldots,X_{10}),$ where $X_{j}=(X_{5(j-1)+1},\ldots
,X_{5(j-1)+5})$, $1\le j \le5$,
and $X_{j}=(X_{3(j-6)+26},\ldots,X_{3(j-6)+28})$, $6 \le j \le10$.
In order to generate $X$, we first simulate 40 random variables
$R_{1},\ldots,R_{40}$, independently from $N(0,1)$. Then, $Z_{j}$,
$j=1,\ldots, 10$ are simulated with a
normal distribution with mean zero and
$\operatorname{cov}(Z_{j_{1}},Z_{j_{2}})=0.6^{|j_{1}-j_{2}|}$. The covariates
$X_{1},\ldots,X_{40}$ are generated as
\begin{eqnarray*}
X_{5(j-1)+k}&=&\frac{Z_{j}+R_{5(j-1)+k}}{\sqrt{2}} , \qquad 1 \leq
j \leq5 ,  1\leq k \leq
5 ,
\\
X_{3(j-6)+25+k}&=&\frac{Z_{j}+R_{3(j-6)+25+k}}{\sqrt{2}} , \qquad 6
\leq j \leq10 ,  1\leq k \leq
3 .
\end{eqnarray*}
The random error $\varepsilon\sim N(0, 3^2)$. The response variable
$Y$ is generated from $Y=\sum_{k=1}^{10}X_{k}\beta_{k}+\varepsilon$, where
$\beta_{1}= (0.5,1,1.5,2,2.5),$   $\beta_{2}=(2,0,0,2,2),$
$\beta_{3} =\cdots=\beta_{5}=(0,0,0,0,0),$
$\beta_{6}=(-1,-2,-3),$
$\beta_{7}=\cdots=\beta_{10}=(0,0,0)$.
\end{example}

\begin{example} \label{e5} In this example, the number of groups is
$p=10$ %is bigger than the sample size $n=200$
and the group sizes differ across groups.
The data are generated the same way as in Example \ref{e4}. However, the
regression coefficients
$\beta_{1} =(0.5,1,1.5,2,2.5),$   $\beta_{2}=(2,2,2,2,2),$
$\beta_{3} =(-1,0,1,2,3),$   $\beta_{4}= (-1.5,2,0,0,0),$
$\beta_{5} =(0,0,0,0,0),$
$\beta_{6} = (2,-2,1),$   $\beta_{7}=(0,-3,1.5),$
$\beta_{8} =(-1.5,1.5,2),$   $\beta_{9}=(-2,-2,-2),$
$\beta_{10}=(0,0,0)$.
\end{example}

\begin{example} \label{e6} In this example, the number of groups
$p=210$ %is bigger than the sample size $n=200$
and the group sizes differ across groups.
The data are generated the same way as in Example \ref{e4}. However, the
regression coefficients
$\beta_{1} =(0.5,1,1.5,2,2.5),$   $\beta_{2}=(2,2,2,2,2),$
$\beta_{3} =(-1,0,1,2,3),$ $\beta_{4}= (-1.5,2,0,0,0),$
$\beta_{5}=\cdots=\beta_{100}=(0,0,0,0,0),$
$\beta_{101} = (2,-2,1),$   $\beta_{102}=(0,-3,1.5),$
$\beta_{103} =(-1.5,1.5,2),$   $\beta_{104}=(-2,-2,-2),$
$\beta_{105} =\cdots=\beta_{210}=(0,0,0)$.
\end{example}

%t1 ###
\begin{table}[b]
\tabcolsep=0pt
\caption{Simulation study by the group Lasso and adaptive group
Lasso for Examples \protect\ref{e1}--\protect\ref{e6}. The true numbers of groups are included in
[] in the first column}\label{tab1}
\begin{tabular*}{\tablewidth}{@{\extracolsep{4in minus 4in}}lllllllllll@{}}
\hline
&\multicolumn{5}{l}{{Group Lasso}} &
\multicolumn{5}{l@{}}{{Adaptive group Lasso}}
\\[-6pt]
&\multicolumn{5}{c}{\hrulefill} &
\multicolumn{5}{c@{}}{\hrulefill}
\\
$\sigma=3$ & mean & med&
ME & \% incl& \% sel & mean & med &  ME & \% incl & \% sel \\
\hline
Ex. 1, [2]& 2.04&2&8.79&100& 96.5& 2.01&2& 8.54&100& 99.5\\
& $(0.18)$&$(2,2)$&$(0.94)$&$(0)$&$(0.18)$&$(0.07)$&$(2,2)$& $(0.90)$&$(0)$&$(0.07)$\\
Ex. 2, [4]& 4.11&4&8.52& 99.5&88.5&4.00&4&8.10&99.5&98.00\\
& $(0.34)$&$(4,4)$& $(0.94)$&$(0.07)$&$(0.32)$&$(0.14)$&$(4,4)$&$(0.87)$&$(0.07)$&$(0.14)$\\
Ex. 3, [4]& 4.00&4& 9.48&93.0& 86.5&3.94&4&8.19&93.0&92.5\\
& $(0.38)$&$(4,4)$& $(1.19)$&$(0.26)$&$(0.34)$&$(0.27)$&$(4,4)$&$(0.96)$&$(0.26)$&$(0.26)$\\
Ex. 4, [3]& 3.17&3& 8.78& 100&85.3& 3.00&3&8.36&100&100\\
& $(0.45)$&$(3,3)$& $(1.00)$& $(0)$&$(0.35)$&$(0)$&$(3,3)$&$(0.90)$&$(0)$&$(0)$\\
Ex. 5, [8]& 8.88 &9& 7.68&100&40.0& 8.03 & 8& 7.58& 100&97.5\\
& $(0.81)$& $(8,10)$& $(0.94)$&$(0)$&$(0.49)$&$(0.16)$&$(8,8)$& $(0.86)$&$(0)$&$(0.16)$\\
Ex 6, [8]& 12.90 &9& 14.61&66.5&7.0&11.49& 8&9.28&66.5&47.0\\
& $(12.42)$&$(8,11)$&$(7.21)$&
$(0.47)$&$(0.26)$&$(12.68)$&$(7,8)$&$(5.79)$&$(0.47)$&$(0.50)$\\
\hline
\end{tabular*}
\end{table}

The results are given in Table \ref{tab1}, based on 400 replications. The
columns in the table include the average number of groups selected
with standard error in parentheses, the median number (`med') of
groups selected with the $25\%$ and $75\%$ quantiles of the number
of selected groups in parentheses, model error (`ME'), percentage of
occasion on which correct groups are included in the selected model
(`$\%$ incl') and percentage of occasions on which the exactly correct
groups are selected (`$\%$ sel'), with standard error in parentheses.

Several observations can be made from Table \ref{tab1}. First, in all six
examples, the adaptive group Lasso performs
better than the group Lasso in terms of model error and the percentage
of correctly selected models. The group Lasso which gives the initial
estimator for the adaptive group Lasso
includes the correct groups with high probability. And the
improvement is considerable for models with different group sizes.
Second, the results from models with
equal group sizes (Examples \ref{e1}, \ref{e2} and \ref{e3}) are better than those from
models with different group sizes (Examples
\ref{e4}, \ref{e5} and \ref{e6}). Finally, when the dimension of the model increases, the
performance of both methods becomes worse.
This is to be expected since selection in models with a larger number
of groups is more difficult.

%s5 ###
\section{Concluding remarks}\label{sec5}

We have studied the asymptotic selection and estimation properties of
the group Lasso and adaptive group Lasso in `large $p$, small $n$'
linear regression models.
For the adaptive group Lasso to be selection consistent,
the initial estimator should possess two properties: (a) it does not
miss important groups and variables; (b) it is estimation consistent,
although it
may not be group-selection or variable-selection consistent.
Under the conditions stated in Theorem \ref{ThmA}, the group Lasso
is shown to satisfy these two requirements.
Thus, the iterated group Lasso procedure, which uses the group Lasso to
achieve dimension reduction and generate the initial estimates and then
uses the adaptive group Lasso
to achieve selection consistency, is an appealing approach to group
selection in high-dimensional settings.

%s6 ###
\section{Proofs}\label{sec6}

We first introduce some notation which will be used in proofs. Let
$\{k\dvt \|\hat{\beta}_{k}\|_{2}>0,k \leq p\}\subseteq A_{1} \subseteq
\{k\dvt
X_{k}^{\prime}(Y-X\hat{\beta})=\lambda\sqrt{d_{k}}\hat{\beta}_{k}/\|
\hat{\beta}_{k}\|_{2}\} \cup
\{1,\ldots, q\}$.
Set $A_{2}=\{1,\ldots,p\}\setminus A_{1}$, $A_{3}=A_{1}\setminus
A_{0}$, $A_{4}=A_{1}\cap A_{0}$, $A_{5}=A_{2}\setminus A_{0}$ and
$A_{6}=A_{2}\cap A_{0}$. Thus, we have $A_{1}=A_{3} \cup A_{4}$,
$A_{3}\cap A_{4}=\varnothing$, $A_{2}=A_{5}\cup A_{6}$ and $A_{5}\cap
A_{6}=\varnothing$. Let $|A_{i}|=\sum_{k\in A_{i}}d_{k}$, $N(A_{i})=
\#\{k\dvt k \in A_{i}\}$, $i=1,\ldots, 6$ and $q_{1}=N(A_{1})$.

\begin{pf*}{Proof of Theorem \ref{ThmA}} The basic idea used
in this proof follows the proof of the rate consistency of the Lasso
in Zhang and Huang (\citeyear{ZH2008}). However, there are many differences in
technical details, for example, in the characterization of the solution via
the Karush--Kuhn--Tucker (KKT) conditions, in the constraint needed for
the penalty level
and in the use of maximal inequalities.

The proof consists of three steps. Step 1 proves some inequalities
related to $q_{1}$, $\tilde{\omega}$ and $\zeta_{2}$. Step~2
translates the results of Step 1 into upper bounds for $\hat{q}$,
$\tilde{\omega}$ and $\zeta_{2}$. Step 3 completes the proof by
showing the probability of the event in Step 2 converging to $1$.
The details of the complete proof are available from the website
\href{http://www.stat.uiowa.edu/techrep}{www.stat.uiowa.edu/techrep}. We will sketch the proof in the
following.

If ${\hat{\beta}}$ is a solution of (\ref{eq:GL1}), then, by the KKT
condition, $X_{k}^{\prime}(Y-X\hat{\beta}) =
\lambda\sqrt{d_{k}}\hat{\beta}_{k}/\|\hat{\beta}_{k}\|_{2}$
$\forall\|\hat{\beta}_{k}\|_{2}>0$ and
$-\lambda\sqrt{d_{k}} \leq X_{k}^{\prime}(Y-X\hat{\beta}) \leq
\lambda\sqrt{d_{k}}$ $\forall
\|\hat{\beta}_{k}\|_{2}=0$. We then have
%
%e6.2 ###
%e6.1 ###
\begin{eqnarray}
\label{eq:c1} & \Sigma_{11}^{-1}S_{A_{1}}/n =
(\beta_{A_{1}}-\hat{\beta}_{A_{1}})+\Sigma_{11}^{-1}\Sigma
_{12}\beta_{A_{2}}
+\Sigma_{11}^{-1}X_{A_{1}}^{\prime}\varepsilon/n ,&
\\
\label{eq:c2} &
n\Sigma_{22}\beta_{A_{2}}-n\Sigma_{21}\Sigma_{11}^{-1}\Sigma
_{12}\beta_{A_{2}}
\leq
C_{A_{2}}-X_{A_{2}}^{\prime}\varepsilon-\Sigma_{21}\Sigma_{11}^{-1}S_{A_{1}}
+\Sigma_{21}\Sigma_{11}^{-1}X_{A_{1}}^{\prime}\varepsilon ,&
\end{eqnarray}
where $S_{A_{i}}=(S_{k_{1}}^{\prime},\ldots,S_{k_{q_{i}}}^{\prime})^{\prime}$,
$S_{k_{i}}=\lambda\sqrt{d_{k_{i}}}s_{k_{i}}$,
$s_{k}=X_{k}^{\prime}(Y-X\hat{\beta})/(\lambda\sqrt{d_{k}})$,
$C_{A_{i}}=(C_{k_{1}}^{\prime},\ldots,\break C_{k_{q_{i}}}^{\prime})^{\prime}$,
$C_{k_{i}}=\lambda
\sqrt{d_{k_{i}}}I(\|\hat{\beta}_{k_{i}}\|_{2}=0)e_{d_{k_{i}}\times
1}$, all the elements of matrix $e_{d_{k_{i}}\times1}$ equal 1,
$k_{i}\in A_{i}$ and $\Sigma_{ij}=X_{A_{i}}^{\prime}X_{A_{j}}/n$.\vspace*{1.5pt}

\textit{Step} 1. Define
\begin{eqnarray*}
V_{1j}=\Sigma_{11}^{-1/2}Q_{A_{j}1}^{\prime}S_{A_{j}}/\sqrt{n} ,\qquad
j=1,3,4 , \qquad\omega_{k}=(I-P_{A_{1}})X_{A_{k}}\beta
_{A_{k}} ,\qquad  k=2,\ldots,6 ,
\end{eqnarray*}
where $Q_{A_{k}j}$ is the matrix representing the selection of
variables in $A_{k}$ from $A_{j}$. Define
$u=X_{A_{1}}\Sigma_{11}^{-1}Q_{A_{4}1}^{\prime}S_{A_{4}}/n-\omega_{2}/
\|X_{A_{1}}\Sigma_{11}^{-1}
Q_{A_{4}1}^{\prime}S_{A_{4}}/n-\omega_{2}\|_{2}$. From $\eqref{eq:c1}$ and
$(\ref{eq:c2})$, we have $V_{14}^{\prime}(V_{13}+V_{14})\leq
S_{A_{4}}^{\prime}Q_{A_{4}1}\Sigma_{11}^{-1}\Sigma_{12}\beta_{A_{2}}+
S_{A_{4}}^{\prime}Q_{A_{4}1}\Sigma_{11}^{-1}X_{A_{1}}^{\prime}\varepsilon/n
+\sqrt{d_{a}}\lambda\sum_{k \in A_{4}}\|\beta_{k}\|_{2}$ and $\|
\omega_{2}\|_{2}^{2} \leq
\beta_{A_{2}}^{\prime}(C_{A_{2}}-X_{A_{2}}^{\prime}\varepsilon-\Sigma
_{21}\Sigma_{11}^{-1}S_{A_{1}}
+\Sigma_{21}\Sigma_{11}^{-1}X_{A_{1}}^{\prime}\varepsilon)$. Then, under
GSC,
\begin{eqnarray}\label{eq:vw}
\|V_{14}\|_{2}^{2}+\|\omega_{2}\|_{2}^{2} &\leq&
(\|V_{14}\|_{2}^{2}+\|\omega_{2}\|_{2}^{2})^{1/2}u^{\prime}\varepsilon
+(\|V_{14}\|_{2}+
\|P_{1}X_{A_{2}}\beta_{A_{2}}\|_{2})\biggl(\frac{\lambda
^{2}d_{a}N(A_{3})}{nc_{*}(|A_{1}|)}\biggr)^{1/2}\nonumber
\\[-8pt]\\[-8pt]
&&{}+\sqrt{d_{a}}\lambda\eta_{1}+\lambda
\sqrt{d_{a}}\|\beta_{A_{5}}\|_{2} .\nonumber
\end{eqnarray}

\textit{Step} 2. Define $B_{1}^2=\lambda^{2}d_{b}q/(nc^{*}(|A_{1}|))$
and $B_{2}^2=\lambda^{2}d_{b}q/(nc_{*}(|A_{0}|\vee
|A_{1}|))$. In this step, we consider the event $|u^{\prime}\varepsilon
|^{2}\leq
(|A_{1}|\vee d_{b})B_{1}^{2}/(4qd_{a})$.
Suppose that the set $A_{1}$ contains all large $\beta_{k} \neq0$. From
$\eqref{eq:vw}$, $\|V_{14}\|_{2}^{2}\leq
B_{1}^{2}+4\sqrt{d_{a}}\lambda
\eta_{1}+4\sqrt{d}\eta_{2}B_{2}+4dB_{2}^{2}$, so we have
%
%e6.3 ###
\begin{equation}\label{eq:Cq}
(q_{1}-q)^{+}\leq
q+\frac{nc^{*}(|A_{1}|)}{\lambda^{2}d_{b}}\Biggl(4\sqrt
{d_{a}}\lambda
\eta_{1}+4\sqrt{\frac{\lambda^{2}d_{a}q} {nc_{*}(|A_{1}|)}}\eta_{2}
+\frac{4\lambda^{2}d_{a}q}{nc_{*}(|A_{1}|)}\Biggr) .
\end{equation}

For general $A_{1}$, let $C_{5}=c^{*}(|A_{5}|)/c_{*}(|A_{1}|\cup
|A_{5}|)$. From $\eqref{eq:vw}$,
%
%e6.4 ###
\begin{equation}\label{eq:Cw}
\|\omega_{2}\|_{2}^{2}\leq
\frac{4}{3}\biggl(\frac{B_{1}^{2}}{2}+dB_{2}^{2}+\sqrt{d}\bigl(1+\sqrt
{C_{5}}\bigr)\eta_{2}B_{2}+2\sqrt{d_{a}}\eta_{1}\biggr)+
\frac{32}{9}dC_{5}B_{2}^{2} .
\end{equation}

From Zhang and Huang (\citeyear{ZH2008}), $\|\omega_{2}\|_{2}^{2} \geq
(\|\beta_{A_{5}}\|_{2}(nc_{*,5})^{1/2}-\eta_{2})^{2}$ and
$\|X_{A_{2}}\beta_{A_{2}}\|_{2}\leq
\eta_{2}+\|X_{A_{5}}\beta_{A_{5}}\|_{2} \leq
\eta_{2}+(nc^{*}(|A_{5}|))^{1/2}\|\beta_{A_{5}}\|_{2}$. By the
Cauchy--Schwarz inequality,
then, we have
\begin{eqnarray}\label{eq:Cr}
 \|\beta_{A_{5}}\|_{2}^{2}nc_{*,5} &\leq&
\Biggl[\frac{4}{3}\lambda\sqrt{\frac{d_{a}q}{nc_{*,5}}}
\biggl(1+\frac{c^{*}(|A_{5}|)}{c_{*}(|A_{1}|)}\biggr)^{1/2}
+2\eta_{2}\Biggr]^{2}\nonumber
\\[-8pt]\\[-8pt]
&&{}+\frac{8}{3}\biggl[\frac{B_{1}^{2}}{4}+\sqrt{d_{a}}\lambda\eta
_{1}+\eta_{2}
\biggl(\frac{\lambda^{2}d_{a}q}{nc_{*}(|A_{1}|)}\biggr)^{1/2}
+\frac{\lambda^{2}d_{a}q}{2nc_{*}(|A_{1}|)}-\frac{3}{4}\eta
_{2}^{2}\biggr] ,\nonumber
\end{eqnarray}
where $c_{*,5}=c_{*}(|A_{1}\cup A_{5}|)$.

\textit{Step} 3. Letting $c_{*}(|A_{m}|)=c_{*}$,
$c^{*}(|A_{m}|)=c^{*}$ for $N(A_{m})\leq q^{*}$,
we have
%
%e6.5 ###
\begin{equation}\label{eq:Cue}
q_{1}\leq N(A_{1}\cup A_{5})\leq q^{*}, \qquad  |u^{\prime}\varepsilon|^{2}\leq
\frac{(|A_{1}|\vee
d_{b})\lambda^{2}d_{b}}{4d_{a}nc^{*}(|A_{1}|)} .
\end{equation}
We have $\bar{c}=C_{5}=c^{*}(|A_{5}|)/c_{*}(|A_{1}|\vee
|A_{5}|)=c^{*}/c_{*}$ and $c_{*,5}=c_{*}(|A_{1}\cup A_{5}|)=c_{*}$.
From $\eqref{eq:Cq}$, $\eqref{eq:Cw}$ and $\eqref{eq:Cr}$,
$(q_{1}-q)^{+}+q \leq M_{1}q$, $\|\omega_{2}\|_{2}^{2} \leq
M_{2}B_{1}^{2}$, $nc_{*}\|\tilde{\gamma}_{A_{5}}\|_{2}^{2}\leq
M_{3}B_{1}^{2}$ when $\eqref{eq:Crk2GL}$ is satisfied. Define
%
%e6.6 ###
\begin{equation}\label{eq:xmm}
x_{m}^{*}\equiv\max_{|A|=m}  \max_{\|U_{A_{k}}\|_{2}=1,
k=1,\ldots,m}\biggl|\varepsilon^{\prime}\frac{X_{A}(X_{A}^{\prime}X_{A})^{-1}
\bar{S}_{A}-(I-P_{A})X\beta}{\|X_{A}(X_{A}^{\prime}X_{A})^{-1} \bar
{S}_{A}-(I-P_{A})X\beta\|_{2}}\biggr|
\end{equation}
for $|A|=q_{1}=m\geq0$,
$\bar{S}_{A}=(\bar{S}_{A_{1}}^{\prime},\ldots,\bar{S}_{A_{m}}^{\prime})^{\prime}$,
where $\bar{S}_{A_{k}}=\lambda\sqrt{d_{A_{k}}}U_{A_{k}}$,
$\|U_{A_{k}}\|_{2}=1$. Let $Q_{A}=X_{A}^{*}(X_{A}^{\prime}X_{A})^{-1}$,
where $X_{k}^{*}=\lambda\sqrt{d_{k}}X_{k}$ for $k \in A$. For a
given $A$, let $V_{lj}=(0,\ldots,0,1,0,\ldots,0)$ be the $|A| \times
1$ vector with the $j$th element in the $l$th group being $1$. Then, by
$\eqref{eq:xmm}$,
\begin{eqnarray*}
 x_{m}^{*} \leq\max_{|A|=m}
\max_{l,j}\biggl\{\biggl|\varepsilon^{\prime}\frac{Q_{A}V_{lj}}{\|
Q_{A}V_{lj}\|_{2}}\biggl|\frac{\|Q_{A}V_{lj}\|_{2}\sum_{l \in
A}\sqrt{d_{l}}}{\|Q_{A}U_{A}\|_{2}}+\biggr|\frac{\varepsilon
^{\prime}(I-P_{A})X\beta}{\|(I-P_{A})X\beta\|_{2}}\biggr|\biggr\} .
\end{eqnarray*}

If we define
$\Omega_{m_{0}}=\{(U,\varepsilon)\dvt x_{m}^{*}\leq
\sigma\sqrt{8(1+c_{0})V^{2}((md_{b}) \vee d_{b})\log(N_{d}\vee
a_{n})}\  \forall m\geq m_{0}\}$, then $(X,\varepsilon)\in\Omega
_{m_{0}}\Rightarrow|u^{\prime}\varepsilon|^{2}\leq(x_{m}^{*})^{2}\leq
(|A_{1}|\vee
d_{b})\lambda^{2}d_{b}/(4d_{a}nc^{*})$ for $N(A_{1})\geq m_{0}\geq0$.
By the definition of $x_{m}^{*}$, it is less than the maximum of
$\left({p\atop m}\right) \sum_{k \in A}d_{k}$ normal variables with mean 0 and
variance $\sigma^{2}V_{\varepsilon}^{2}$, plus the maximum of
$\left({p\atop m}\right)$ normal variables with mean 0 and variance
$\sigma^{2}$. It follows that
$P\{(X,\varepsilon)\in\Omega_{m_{0}}\} \rightarrow1$
when $\eqref{eq:Cue}$ holds. This completes the sketch of the proof
of Theorem \ref{ThmA}.\end{pf*}

\begin{pf*}{Proof of Theorem \ref{ThmB}} Consider the case
when $\{c^{*}, c_{*}, r_{1}, r_{2}, c_{0},
\sigma, d\}$ are fixed. The required configurations in Theorem \ref
{ThmA} then become
%
%e6.7 ###
\begin{equation}
\label{Req}
M_{1}q+1<q^{*} ,\qquad\eta_{1}\leq\frac{r_{1}^{2}}{c^{*}}\frac
{q\lambda}{n} , \qquad\eta_{2}^{2}\leq
\frac{r_{2}^{2}}{c^{*}}\frac{q\lambda^{2}}{n} .
\end{equation}

Let $A_{1}=\{k\dvt \|\hat{\beta}_{k}\|_{2}>0   \mbox{ or }   k \notin
A_{0}\}$. Define
$v_{1}=X_{A_{1}}(\hat{\beta}_{A_{1}}-\beta_{A_{1}})$ and\vspace*{-1pt}
$g_{1}=X_{A_{1}}^{\prime}(Y-X\hat{\beta})$. We then have $\|v_{1}\|
_{2}^{2}\geq
c_{*}n\|\hat{\beta}_{A_{1}}-\beta_{A_{1}}\|_{2}^{2}$,
$(\hat{\beta}_{A_{1}}-\beta_{A_{1}})^{\prime}g_{1}
=v_{1}^{\prime}(X\beta-X_{A_{1}}\beta_{A_{1}}+\varepsilon)-\|v_{1}\|_{2}^{2}$
and $\|g_{1}\|_{\infty}\leq
\max_{k,\|\hat{\beta}_{k}\|_{2}>0}\|\lambda\sqrt{d_{k}}\hat{\beta
}_{k}/\|\hat{\beta}_{k}\|_{2}\|_{\infty}=\lambda
d_{a}$. Therefore, $\|v_{1}\|_{2}\leq
\eta_{2}+\|P_{A_{1}}\varepsilon\|_{2} +\lambda
\sqrt{d_{a}N(A_{1})/(nc_{*})}$. Since
$\|P_{A_{1}}\varepsilon\|_{2}\leq2\sigma\sqrt{N(A_{1})\log(N_{d})}$
with probability converging to $1$ under the normality assumption,
$\|X(\hat{\beta}-\beta)\|_{2}\leq
2\eta_{2}+\|P_{A_{1}}\varepsilon\|_{2}
+\lambda\sqrt{d_{a}N(A_{1})/(nc_{*})}$. We then have
%
%e6.8 ###
\begin{equation}
\label{eq:A1}
\biggl(\sum_{k \in A_{1}}\|\hat{\beta}_{k}-\beta_{k}\|_{2}^{2}\biggr)^{1/2}\leq
\frac{\|v_{1}\|_{2}}{\sqrt{nc_{*}}} \leq
\frac{1}{\sqrt{nc_{*}}}\bigl(\eta_{2}+2\sigma\sqrt{N(A_{1})\log
(N_{d})}+\sqrt{dM_{1}\bar{c}}B_{1}\bigr) .
\end{equation}
Since $A_{2} \subset A_{0}$, by the second inequality in
(\ref{Req}), $\# \{k \in A_{0}\dvt \|\beta_{k}\|_{2}>\lambda/n\} \leq
r_{1}^{2}q/c^{*}\sim \mathrm{O}(q)$. By the SRC and the third inequality
in (\ref{Req}), $\sum_{k \in A_{0}}\|\beta_{k}\|_{2}^{2}I \{\|\beta
_{k}\|_{2}>\lambda/n\} \leq\sum_{k \in
A_{0}}\|X_{k}\beta_{k}\times I\{\|\beta_{k}\|_{2}>\lambda/n\}\|
_{2}^{2}/(nc_{*}) \leq r_{2}^{2}q\lambda^{2}/(n^{2}c_{*}c^{*})$ and
$\sum_{k \in A_{0}}\|\beta_{k}\|_{2}^{2}I\{\|\beta_{k}\|_{2}\leq
\lambda/n\} \leq
r_{1}^{2}q\lambda^{2}/(c^{*}n^{2})$.
From (\ref{eq:A1}), we then have
\begin{eqnarray*}
\|\hat{\beta}-\beta\|_{2} &\leq&\frac{1}{\sqrt{nc_{*}}}
\bigl(2\sigma\sqrt{M_{1}\log(N_{d})q}+\bigl(r_{2}
+\sqrt{dM_{1}\bar{c}}\bigr)B_{1}\bigr)+ \sqrt{\frac
{c_{*}r_{1}^{2}+r_{2}^{2}}{c_{*}c^{*}}}\frac{\sqrt{q}\lambda}{n} ,
\\
\|X\hat{\beta}-X\beta\|_{2}&\leq&2\sigma\sqrt{M_{1}\log
(N_{d})q}+\bigl(2r_{2}+\sqrt{dM_{1}\bar{c}}\bigr)B_{1} .
\end{eqnarray*}
This completes the proof of Theorem \ref{ThmB}.
\end{pf*}

\begin{pf*}{Proof of Theorem \ref{ThmC}} Let
$\hat{u}={\hat{\beta}}-\beta$, $W=X^{\prime}\varepsilon/\sqrt{n}$,
$V(u)=\sum_{i=1}^{n}[(\varepsilon_{i}-x_{i}u)^{2}-\varepsilon_{i}^{2})]+
\sum_{k=1}^{p}\lambda_{k}\sqrt{d_{k}}\|u_{k}+\beta_{k}\|_{2}$ and
$\hat{u}=\min_{u}(\varepsilon-Xu)^{\prime}(\varepsilon-Xu)+
\sum_{k=1}^{p}\lambda_{k}\sqrt{d_{k}}\|u_{k}+\beta_{k}\|_{2}$, where
$\lambda_k = \tilde{\lambda}/\|\tilde{\beta}_k\|_2$. By the KKT conditions,
if there exists $\hat{u}$ such that
%
%e6.10 ###
%e6.9 ###
\begin{eqnarray}
\label{eq:o1}
&\Sigma_{A_{0}^{c}A_{0}^{c}}\bigl(\sqrt{n}\hat{u}_{A_{0}^{c}}\bigr)-W_{A_{0}^{c}}
= -S_{A_{0}^{c}}/\sqrt{n} , \qquad\|\hat{u}_{k}\|_{2}\leq
\|\beta_{k}\|_{2}\qquad\mbox{for }k \in A_{0}^{c} ,&
\\
\label{eq:o3}
& -C_{A_{0}}/\sqrt{n} \leq
\Sigma_{A_{0}A_{0}^{c}}\bigl(\sqrt{n}\hat{u}_{A_{0}^{c}}\bigr)-W_{A_{0}} \leq
C_{A_{0}}/\sqrt{n} ,&
\end{eqnarray}
then $\|\hat{\beta}_{k}\|_{2}\neq0$ for $k= 1, \ldots, q$ and
$\|\hat{\beta}_{k}\|_{2}=0$ for $k=q+1, \ldots, p$.\vspace*{1pt}

From $\eqref{eq:o1}$ and $\eqref{eq:o3}$,
$(\sqrt{n}\hat{u}_{A_{0}^{c}})-\Sigma_{A_{0}^{c}A_{0}^{c}}^{-1}W_{A_{0}^{c}}
= -\frac{1}{\sqrt{n}}\Sigma_{A_{0}^{c}A_{0}^{c}}^{-1}
S_{A_{0}^{c}}$ and
$\Sigma_{A_{0}A_{0}^{c}}(\sqrt{n}\hat
{u}_{A_{0}^{c}})-W_{A_{0}}=-n^{-1/2}X_{A_{0}}^{\prime}(I-P_{A_{0}^{c}})\varepsilon
-n^{-1/2}\Sigma_{A_{0}A_{0}^{c}}
\Sigma_{A_{0}^{c}A_{0}^{c}}^{-1}S_{A_{0}^{c}}$. Define the events
\begin{eqnarray*}
E_1&=&\bigl\{
n^{-1/2}\|(\Sigma
_{A_{0}^{c}A_{0}^{c}}^{-1}{X_{A_{0}^{c}}^{\prime}\varepsilon})_{k}\|_{2}<
\sqrt{n}\|\beta_{k}\|_{2}-
n^{-1/2}\|(\Sigma_{A_{0}^{c}A_{0}^{c}}^{-1}S_{A_{0}^{c}})_{k}\|_{2},
  k \in A_{0}^{c}\bigr\} ,
  \\
E_2&=&\bigl\{
n^{-/2}\bigl\|\bigl(X_{A_{0}}^{\prime}(I-P_{A_{0}^{c}})\varepsilon\bigr)_{k}\bigr\|_{2}<
n^{-1/2}\|C_{k}\|_{2}-n^{-1/2}
\|(\Sigma_{A_{0}A_{0}^{c}}\Sigma
_{A_{0}^{c}A_{0}^{c}}^{-1}S_{A_{0}^{c}})_{k}\|_{2},
  k \in A_{0}\bigr\} ,
\end{eqnarray*}
where $(\cdot)_{k}$ denotes the $d_{k}$-dimensional subvector of
the vector $(\cdot)$ corresponding to the $k$th group. We then have
$P(\|\hat{\beta}_{k}\|_{2}\neq0$, $k\in A_0, \mbox{ and }
\|\hat{\beta}_{k}\|_{2}=0, k\notin A_0)\geq P(E_1 \cap E_2)$ and
$P(E_1 \cap E_2)=1-P(E_1^{c}\cup E_2^{c})\ge
1-P(E_1^{c})-P(E_2^{c})$.

First, we consider $P(E_1^{c})$. Define $R=\{\|\tilde{\beta}_{k}\|_{2}^{-1}
\leq c_{1}\theta_{b}^{-1}, k \in A_{0}^{c} \}$, where $c_{1}$ is a
constant. $P(E_1^{c})=P(E_1^{c} \cap R)+P(E_1^{c}\cap R^{c})
\leq P(E_1^{c} \cap R)+P(R^{c})$. By (C2),\vspace*{1pt} $P(R^{c})\rightarrow
0$. Let $N_{q}=\sum_{k=1}^{q}d_{k}$, $\tau_{1}\leq\cdots\leq
\tau_{N_{q}}$ be the eigenvalues of $\Sigma_{A_{0}^{c}A_{0}^{c}}$
and $\gamma_{1},\ldots, \gamma_{N_{q}}$ be the associated
eigenvectors. The $j$th
element in the $l$th group of vector
$\Sigma_{A_{0}^{c}A_{0}^{c}}^{-1}S_{A_{0}^{c}}$ is
$u_{lj}=\sum_{l^{\prime}=1}^{N_{q}}\tau_{l^{\prime}}^{-1}
(\gamma_{l^{\prime}}^{\prime}S_{A_{0}^{c}})\gamma_{lj}$. By the
Cauchy--Schwarz inequality, $u_{lj}^{2} \leq
\tau_{1}^{-2}\sum_{l=1}^{N_{q}}\|\gamma_{l}\|_{2}^{2}\|
S_{A_{0}^{c}}\|_{2}^{2}
= \tau_{1}^{-2}N_{q}\|S_{A_{0}^{c}}\|_{2}^{2} \leq
\tau_{1}^{-2}N_{q}(\sum_{k=1}^{q}\lambda_{k}^{2}d_{k})$. Therefore,
$\|u_{k}\|_{2}^{2}\leq
d_{k}\tau_{1}^{-2}q^{2}d_{a}^{2}(\tilde{\lambda}
c_{1}\theta_{b}^{-1})^{2}$.

If we define $\upsilon_{A_{0}^{c}}=
\sqrt{n}\theta_{b}-n^{-1/2}c_{1}\tau_{1}^{-1}qd_{a}^{3/2}\tilde
{\lambda}
\theta_{b}^{-1}$,
$\eta_{A_{0}^{c}}=n^{-1/2}\Sigma
_{A_{0}^{c}A_{0}^{c}}^{-1}X_{A_{0}^{c}}^{\prime}\varepsilon$,
$\xi_{A_{0}}=n^{-1/2}\times X_{A_{0}}^{\prime}(I-P_{A_{0}^{c}})\varepsilon$,
$C_{A_{0}^{c}}=\{\max_{k \in A_{0}^{c}}\|\eta_{k}\|_{2} \geq
\upsilon_{A_{0}^{c}}\}$, then $P(E_1^{c})\leq P(C_{A_{0}^{c}})$. By
Lemmas 1 and~2 of Huang, Ma and Zhang (\citeyear{HHM2008}), $P(C_{A_{0}}^{c}) \leq
K(d_{a}\log q)^{1/2}/\upsilon_{A_{0}^{c}}$, where $K$ is a constant,
$k(d_{a}\log q)^{1/2}/\upsilon_{A_{0}^{c}} \rightarrow0$ from
(C3). We then have $P(E_1^{c}\cap R)\rightarrow0$,
$P(E_1^{c})\rightarrow0.$

Next, we consider $P(E_2^{c})$. Similarly as above, define
$D=\{\|\tilde{\beta}_{k}\|_{2}^{-1} > r_{n}, k \in A_{0}\} \cap R$.
$P(E_2^{c}) \leq
P(E_2^{c}\cap D)+P(D^{c})$. By (C2), $P(D^{c}) \rightarrow0$.
$|\sum_{l=1}^{N_{q}}\sum_{i=1}^{n}(X_{A_{0}})_{ij}
(X_{A_{0}^{c}})_{il}u_{l}|\leq\sum_{l=1}^{N_{q}}|u_{l}/n| \leq
\tau^{-1}_{1}q^{2}d_{a}^{2}\tilde{\lambda} c_{1}\theta_{b}^{-1}$,
where $u_{l}$ is the $l$th element of vector
$\Sigma_{A_{0}^{c}A_{0}^{c}}^{-1}S_{A_{0}^{c}}$. If we define\vspace*{-2pt}
$\upsilon_{A_{0}}=n^{-1/2}\tilde{\lambda} r_{n}\sqrt{d_{b}}
-n^{-1/2}\tau_{1}^{-1}q^{2}d_{a}^{5/2}\tilde{\lambda}
c_{1}\theta_{b}^{-1}$, $C_{A_{0}}=\{\max_{k \in
A_{0}}\|\xi_{k}\|_{2} > \upsilon_{A_{0}}\}$, then $P(Q^{c})\leq
P(C_{A_{0}})$, $P(C_{A_{0}}) \leq K(d_{a}\log
(p-q))^{1/2}/\upsilon_{A_{0}}$. $K(d_{a}\log(p-q))^{1/2}/\upsilon
_{A_{0}} \rightarrow0$ from (C3).\vspace*{1pt}
We then have $P(E_2^{c}\cap D) \rightarrow0$,
$P(E_2^{c})\rightarrow0.$
This completes the proof of Theorem \ref{ThmC}.
\end{pf*}

\begin{pf*}{Proof of Theorem \ref{ThmD}}
If we let $\hat{A}=\{k\dvt \|\tilde{\beta}_{k}\|_{2}>0, k=1,\ldots, p\}
$, %Then
then $\sum_{k \in{\hat{A}}^{c}}\|\hat{\beta}^{*}_{k}\|_{2}=0$, the
dimension of our problem (\ref{eq:ModelAGL}) is reduced to $\hat{q}$,
$\hat{q}\leq q^{*}$ and $\hat{A}_{c}\subset A_{0}$. By the
definition of $\hat{\beta}^{*}$,
we have
%
%e6.12 ###
%e6.11 ###
\begin{eqnarray}
\label{eq:41} &\displaystyle\frac{1}{2}\|Y-X_{\hat{A}}\hat{\beta}^{*}_{\hat
{A}}\|
^{2}_{2}+\tilde{\lambda}\sum_{k \in
\hat{A}}\frac{\sqrt{d_{k}}}{\|\tilde{\beta}_{k}\|_{2}}\|\hat
{\beta}_{k}^{*}\|
_{2} \leq
\frac{1}{2}\|Y-X_{\hat{A}}\beta_{\hat{A}}\|^{2}_{2}+\tilde{\lambda
}\sum_{k \in
\hat{A}}\frac{\sqrt{d_{k}}}{\|\tilde{\beta}_{k}\|_{2}}\|\beta
_{k}\|_{2} ,\quad &
\\
\label{eq:42} &\displaystyle\eta^{*}=\tilde{\lambda} \sum_{k \in\hat{A}}\frac
{\sqrt{d_{k}}}{\|\tilde{\beta}_{k}\|_{2}}
(\|\beta_{k}\|_{2}-\|\hat{\beta}^{*}_{k}\|_{2})\leq\tilde{\lambda
} \sum
_{k \in
\hat{A}}\frac{\sqrt{d_{k}}}{\|\tilde{\beta}_{k}\|_{2}}\|\hat
{\beta}
_{k}^{*}-\beta_{k}\|_{2} .&
\end{eqnarray}

If we let
$\delta_{\hat{A}}=\Sigma_{\hat{A}\hat{A}}^{1/2}(\hat{\beta
}^{*}_{\hat
{A}}-\beta_{\hat{A}})$ and
$D=\Sigma_{\hat{A}\hat{A}}^{-1/2}X_{\hat{A}}^{\prime}$, then
$\|Y-X_{\hat{A}}\hat{\beta}^{*}_{\hat{A}}\|_{2}^{2}/2-\|Y-X_{\hat
{A}}\beta_{\hat{A}}\|_{2}^{2}/2
=\delta_{\hat{A}}^{\prime}\delta_{\hat{A}}/2-(D\varepsilon)^{\prime}\delta
_{\hat{A}}$.
By (\ref{eq:41}) and (\ref{eq:42}),
$\delta_{\hat{A}}^{\prime}\delta_{\hat{A}}/2-(D\varepsilon)^{\prime}
\delta_{\hat{A}}-\eta^{*}\leq0$, so
$\|\delta_{\hat{A}}-D\varepsilon\|_{2}^{2}-\|D\varepsilon\|
_{2}^{2}-2\eta^{*}
\leq0$. By the triangle inequality, $\|\delta_{\hat{A}}\|_{2}\leq
\|\delta_{\hat{A}}-D\varepsilon\|_{2}+\|D\varepsilon\|_{2}$. Thus,
$\|\delta_{\hat{A}}\|_{2}^{2}\leq
6\|D\varepsilon\|_{2}^{2}+6\eta^{*}$.

Let $D_{i}$ be the $i$th
column of $D$. $ E(\|D\varepsilon\|_{2}^{2})=\sigma^{2}\operatorname{tr}(D^{\prime}D)
=\sigma^{2}\hat{q}$. Then, with probability converging to $1$,
$\|\hat{\beta}_{\hat{A}}-\beta_{\hat{A}}\|^{2}_{2}\leq
6\sigma^{2}M_{1}q/(nc_{*})+(\tilde{\lambda}
\sqrt{d_{a}}/(\xi_{b}\theta_{b}nc_{*}))^{2}/2+
\|\hat{\beta}_{\hat{A}}-\beta_{\hat{A}}\|^{2}_{2}/2$. Thus, for
$\tilde{\lambda}=n^{\alpha}$ for some $0<\alpha<1/2$, with
probability converging to $1$,
\begin{eqnarray*}
\|\hat{\beta}_{\hat{A}}-\beta_{\hat{A}}\|_{2}\leq
\sqrt{\frac{6\sigma^{2}M_{1}}{c_{*}}\frac{q}{n}+\frac{d_{a}}{(\xi
_{b}\theta_{b}c_{*})^{2}}\biggl(\frac{\tilde{\lambda}}{n}\biggr)^{2}}
\sim \mathrm{O}\Biggl(\sqrt{\frac{q}{n}}\Biggr)
\end{eqnarray*}
and
$\|X_{\hat{A}}\hat{\beta}_{\hat{A}}-X_{\hat{A}}\beta_{\hat{A}}\|
_{2}\leq
\sqrt{nc^{*}}\|\hat{\beta}_{\hat{A}}-\beta_{\hat{A}}\|_{2} \sim
\mathrm{O}(\sqrt{q})$. This completes the proof of Theorem~\ref{ThmD}.
\end{pf*}

\section*{Acknowledgements}

 The authors are grateful to Professor Cun-Hui Zhang for
sharing his insights into the problem and related topics. The work of
Jian Huang is supported in part by NIH Grant R01CA120988 and NSF Grants
DMS-07-06108 and 0805670.

\printhistory

\end{document}